\documentclass[10pt,a4paper]{amsart}
\usepackage{amsthm, amsfonts, amssymb,latexsym}
\usepackage{enumerate, color}
\usepackage[latin1]{inputenc}

\usepackage{amsmath}
\usepackage{amsfonts}
\usepackage{amssymb}
\usepackage{graphicx}
\usepackage{fullpage}
\usepackage{bbm}
\usepackage[T1]{fontenc}
\usepackage{parskip}
\usepackage{pdfpages}
\usepackage{url}

\newtheorem{prop}{Proposition}
\newtheorem{lemma}{Lemma}
\newtheorem{theorem}{Theorem}
\newtheorem{cor}{Corollary}
\newtheorem{defi}{Definition}

\numberwithin{cor}{section} \numberwithin{theorem}{section}
\numberwithin{lemma}{section} \numberwithin{prop}{section}

\newcommand{\A}{A}

\newcommand{\EE}{\mathbb{E}}

\title{Ratio sets of random sets}

\author[J. Cilleruelo]{Javier Cilleruelo}
\address{J. Cilleruelo: Instituto de Ciencias Matem\'aticas
(CSIC-UAM-UC3M-UCM) and Departamento de Matem\'aticas, Universidad
Aut\'onoma de Madrid, 28049 Madrid, Spain }
\email{franciscojavier.cilleruelo@uam.es}
\author[J. Guijarro-Ordóñez]{Jorge Guijarro-Ordóñez}
\address{J. Guijarro-Ordóñez: Instituto de Ciencias Matem\'aticas
(CSIC-UAM-UC3M-UCM) and Universidad
Complutense de Madrid, 28040 Madrid, Spain } \email{jorge.guijarro.ord@hotmail.com}

\begin{document}

\begin{abstract}
We study the typical behavior of the size of the ratio set $A/A$ for a random subset
$\A\subset \{1,\dots , n\}$. For example, we prove that $|A/A|\sim \frac{2\text{Li}_2(3/4)}{\pi^2}n^2 $ for almost all subsets $A\subset\{1,\dots ,n\}$. We also prove that the proportion of visible lattice points in the lattice $A_1\times\cdots \times A_d$, where $A_i$ is taken at random in $[1,n]$ with $\mathbb P(m\in A_i)=\alpha_i$ for any $m\in [1,n]$, is asymptotic to a constant $\mu(\alpha_1,\dots,\alpha_d)$ that involves the polylogarithm of order $d$.
\end{abstract}
\maketitle
\section{Introduction}\label{sec:intro}
Given a set of positive integers $A$, we say that a lattice point $P\in \mathbb N\times \mathbb N$ is visible in the lattice $A\times A$ if the line segment connecting $(0,0)$ and $P$ does not contain more lattice points of $A\times A$. We denote by $\text{visible}(A\times A)$ the set of the visible lattice points in $A\times A$ and denote by $A/A$ the ratio set $A/A=\{a/a':\ a,a'\in A\}$. Each visible lattice point in the lattice $A\times A$ may be identified with an element of $A/A$ and then we have that $|\text{visible}(A\times A)|=|A/A|$. It is well known that the set of visible lattice points in the plane has density $6/\pi^2$; so, if we write $I_n=\{1,\dots,n\}$, we have that $|\text{visible}(I_n\times I_n)|=|I_n/I_n|\sim \frac 6{\pi^2}n^2$. 

In the present paper we study the typical size of  $\text{visible}(A\times A)$  for a random set $\A$ in $\{1,\dots ,n\}$ when $n\to \infty$ or, equivalently, the typical size of the ratio set $A/A$. We consider two natural probabilistic models.

In the first one, denoted by $B(n,\alpha)$, each element in $\A$ is chosen independently at random in $\{1,\dots, n\}$ with probability  $\alpha$. Then we have the following
 \begin{theorem}\label{TH1} Let  $\alpha\in (0,1)$ and consider a random subset $A$ in $B(n,\alpha)$.  Then, with probability $1-o(1)$ when $n\rightarrow\infty$, 
$$|A/A|\sim \mu_2(\alpha)n^2,
$$
where $$\mu_2(\alpha)=\frac{\alpha^2}{\zeta(2)}\frac{\mathrm{Li}_2(1-\alpha^2)}{1-\alpha^2} $$ and $\mathrm{Li}_2(z)=\sum_{k\ge 1}\frac{z^k}{k^2}$ is the dilogarithm function.
\end{theorem}
Notice that the case $\alpha=1$ corresponds to take $A=I_n$ and its asymptotic estimate appears as the limiting case, as $\alpha$ tends to 1, in Theorem \ref{TH1}, since $\lim_{\alpha \to 1}\mu_2(\alpha)=\zeta^{-1}(2)=6/\pi^2$. 

Furthermore, when  $\alpha=1/2$, all the subsets $A\subset \{1,\dots,n\}$ are chosen with the same probability and Theorem \ref{TH1} gives the following result.
\begin{cor}We have that
$$\frac 1{2^n}\sum_{A\subset \{1,\dots,n\}}|A/A|\sim \mu_2(1/2)n^2.$$
Furthermore, for almost all sets $\A\subset \{1,\dots,n\}$ we have that $$|A/A|\sim \mu_2(1/2)n^2.$$
\end{cor}
A strong convergence version of Theorem \ref{TH1} is also possible.
\begin{cor}\label{strong}
Let $A$ be a random infinite sequence of positive integers where all the events $m\in A$ are independent and $\mathbb P(m\in A)=\alpha$ for any positive integer $m$. Let $A_n$ be the random set $A\cap [1,n]$. Then we have
$$\mathbb P\left (\lim_{n\to \infty}\frac{|A_n/A_n|}{n^2}= \mu_2(\alpha) \right )=1.$$
\end{cor}

For a given positive integer $k=k(n)$,  typically $k\asymp n$, we consider the second model, where each subset of $k$ elements is chosen uniformly at random among  all sets of  size $k$ in $ \{1,\dots, n\}$. We denote this model by $S(n,k)$.

When $k/n\sim \alpha$, the heuristic suggests that both models, $S(n,k)$ and $B(n,\alpha)$, are quite similar. Indeed, this is the strategy we follow to prove Theorem \ref{TH2}.
\begin{theorem}\label{TH2}Let  $k\asymp n$ and consider a random subset $A$ in $S(n,k)$.  Then, with probability $1-o(1)$ when $n\rightarrow\infty$, 
$$|A/A|\sim \mu_2(k/n)n^2,$$
where $\mu_2$ denotes the function of the same name introduced in Theorem \ref{TH1}. 


\end{theorem}
The case $k=n$, which corresponds to the classical result in which $A=I_n$, is also obtained as a limiting case in Theorem \ref{TH2} in the sense that $\lim_{k/n\to 1}\mu_2(k/n)=\zeta^{-1}(2)=6/\pi^2$.

The next theorem deals with the size of $A/A$ in the special case in which $A$ is an arithmetic progression:
\begin{theorem}\label{TH3}
Let $A$ be the set of integers congruent to $a\pmod q$ in $\{1,\dots,n\}$. We have
$$|A/A|\sim c_q\frac {n^2}{q^2}$$
when $n\rightarrow\infty$, where
$$c_q=\frac 6{\pi^2}\prod_{p\mid q}\left (1-\frac 1{p^2}\right )^{-1}\sum_{\substack{1\le t\le q\\ (t,q)=1}}\frac 1{t^2}.$$
\end{theorem}
Note that, when $q=1$, we recover once more the classical result for the plane. In the opposite case, it is easy to see that $c_q\to 1$ when $q\to \infty$. Moreover, if we let $\alpha=1/q$ in Theorem \ref{TH1}, we realize that arithmetic progressions are an \textit{atypical} set for this problem.

The results above might suggest that we always have that $|A/A|\asymp|A|^2$ when $A\subset\{1,\dots ,n\}$ has positive density; or, equivalently, that $|A/A|\gg|A|^2$ since the bound $|A/A|\leq |A|^2$ is trivial. However, this intuition is false, as the next theorem shows.

\begin{theorem}\label{TH5}
For any $\epsilon>0$, there exists $\alpha>0$ such that, for all sufficiently large $n$, there exists a subset $A$ of the integers in $[1, n]$ satisfying $|A|\ge \alpha n$ and $|A/A|< \epsilon|A|^2$.
\end{theorem}
This was proved in a different setting in \cite{Conjuntos}, but here we give an alternative proof, which is simpler and more compact. It is true, however, that $|A/A|\gg_{\delta}\alpha^{\delta}|A|^2$ for any $\delta>0$ and for any $A\subset\{1,\dots ,n\}$ of density $\alpha>0$, as proved in \cite{Conjuntos}.

Finally, we might ask what happens with visible lattice points in multidimensional spaces, considering  sets of not necessarily the same density in the two models above. The answer to this last question is given by Theorem \ref{TH4}.
 Notice that $|A_1/A_2|=|\text{visible}(A_1\times A_2)|$, but there is no ratio set version for dimension $d\ge 3$. This is the reason why we state Theorem \ref{TH4} in terms of visible lattice points.
\begin{theorem}\label{TH4}
 Let  $\alpha_1,\dots,\alpha_d\in (0,1)$ and consider random subsets $A_i$ in $B(n,\alpha_i),\ i=1, \dots, d$.  Then
 $$|{\rm visible}(A_1\times \cdots \times A_d)|\sim \mu(\alpha_1,\dots,\alpha_d)n^d$$ with probability $1-o(1)$ when $n\to \infty$, 
 where $$\mu(\alpha_1,\dots,\alpha_d)=\frac{\alpha_1\cdots \alpha_d} {\zeta(d)}\frac{{\rm Li}_d(1-\alpha_1\cdots \alpha_d)}{1-\alpha_1\cdots \alpha_d}$$ and $\mathrm{Li}_d(z)=\sum_{k\ge 1}\frac{z^k}{k^d}$ is the polylogarithm of order $d$.
\end{theorem}
In particular, this provides yet another instance in which polylogarithms occur. See \cite{Z} for more applications of the polylogarithms. Moreover, when $\alpha_i=1$ for all $1\leq i\leq d$ we recover the classical result (\cite{Nymann}) that the probability that $d$ positive integers are relatively prime is $1/\zeta(d)$, since $\lim_{\alpha_1\cdots \alpha_d\to 1}\mu(\alpha_1,\dots,\alpha_d)=1/\zeta(d)$. 

Theorem \ref{TH4} also works for the probabilistic models $S(n,k_i),\ i=1,\dots, d$. However, we have decided to omit its proof for the sake of brevity, since the ideas involved are those in the proof of Theorem \ref{TH2}. 

There is also a strong convergence version similar to that of Corollary \ref{strong}: \textit{If $A_1,\dots,A_d$ are infinite random sequences of positive integers such that all the events $m\in A_i$ are independent and $\mathbb P(m\in A_i)=\alpha_i,\ i=1,\dots,d$, the random variables $X_n=|\mathrm{visible}(A_1\times \cdots \times A_d\cap [1,n]^d)|$ satisfy that
$$\mathbb P\left(\lim_{n\to \infty}\frac{X_n}{n^d}=\mu(\alpha_1,\dots,\alpha_d)\right)=1.$$}
The proof is similar to that of Corollary \ref{strong}, so details will be omitted.

One last remark is in order: Theorem \ref{TH4} is not a generalization of Theorem \ref{TH1}. Note that, in Theorem \ref{TH1}, we consider the random variable $|\text{visible}(A\times A)|$ where $A$ is a random set in $S(n,\alpha)$; whereas, in Theorem \ref{TH4}, when $d=2$ and $\alpha_1=\alpha_2=\alpha$ we deal with the random variable $|\text{visible}(A_1\times A_2)|$ where $A_1,A_2$ are random sets in $S(n,\alpha)$. However, the natural generalization of Theorem \ref{TH1} also holds, even in the strong convergence version:

\textit{Let $A$ be an infinite random sequence of positive integers where all the events $m\in A$ are independent and $\mathbb P(m\in A)=\alpha$, and let $X_n=|\mathrm{visible}(A\times \stackrel{d}{\cdots} \times A)\cap [1,n]^d|$.
Then
$$\mathbb P\left(\lim_{n\to \infty}\frac{X_n}{n^d}=\mu_d(\alpha)\right)=1,$$ where $\mu_d(\alpha)=\mu(\alpha,\stackrel{d}{\dots},\alpha)$.}

Again the proof follows the same steps that those of Theorem \ref{TH1} and Corollary \ref{strong}.

We now pass to the proofs of the five previous theorems, organized accordingly in five different sections.

\section{The size of $A/A$ for random sets in $B(n,\alpha)$. Proof of Theorem
\ref{TH1}}\label{sec:expectation}

\subsection{Expectation}
First of all, we give an explicit expression for the expected value of
the random variable $X=|A/A|$, where $A$ is a random set in $B(n,\alpha)$.
\begin{prop}\label{expectation} For the random variable $X=|A/A|$ in $B(n,\alpha)$ we have
$$\EE\left(X\right)=\frac {6}{\pi^2}(n\alpha)^2\frac{{\rm Li}_2(1-\alpha^2)}{1-\alpha^2}+O(n\log^2 n).$$
\end{prop}
\begin{proof}
Linearity of expectation and the symmetry with respect to the line $r=s$ give the equality
\begin{equation}\label{eq1}
\mathbb{E}(X)=\sum_{\substack{r,s\leq n\\
            (r,s)=1}}\mathbb{P}(r/s\in A/A)=2\sum_{\substack{r< s\leq n\\
            (r,s)=1}}\mathbb{P}(r/s\in A/A)+1.
\end{equation}      
Moreover, for $r<s$ with $(r,s)=1$,
$$\mathbb{P}(r/s \notin A/A)=\mathbb{P}\left(\bigcap_{t\leq n/s}\lbrace(rt,st) \notin A\times A\rbrace\right)=\mathbb{P}\left(\bigcap_{t\leq n/s}E^c_t\right),$$
where $E_t$ stands for the event $\lbrace rt, st\in A\rbrace$. 

Clearly, these events are independent if and only if there do not exist $t,t'\leq n/s$ such that $rt'=st$. Since $(r,s)=1$, the former condition implies that $s|t'$, so $s\leq t'$. Thus, the inequality $s>\sqrt{n}$ entails the independence of the events, because otherwise we are led to the contradiction $\sqrt{n}< s\leq t'\leq n/s<\sqrt{n}$. Hence, if $s>\sqrt{n}$,
\begin{equation}\label{eq2}
\mathbb{P}(r/s \notin A/A)=\mathbb{P}\left(\bigcap_{t\leq n/s}E^c_t\right)=\prod_{t\leq n/s}\big (1-\mathbb{P}(rt\in A)\mathbb{P}(st\in A)\big )=(1-\alpha^2)^{[n/s]}.
\end{equation}
We consequently split the sum in ~\eqref{eq1} in two parts:
$$\mathbb{E}(X)=2\sum_{\substack{\sqrt{n}< s\leq n\\
            r<s, (r,s)=1}}\mathbb{P}(r/s\in A/A)+2\sum_{\substack{r< s\leq \sqrt{n}\\
            (r,s)=1}}\mathbb{P}(r/s\in A/A)+1.$$
In the first one, we have independent events and equation ~\eqref{eq2} holds. In the second one, we simply bound the probabilities by 1. Therefore,
\begin{eqnarray*}
\mathbb{E}(X)&=&2\sum_{\substack{\sqrt{n}< s\leq n\\
            r< s,\ (r,s)=1}}\left(1-(1-\alpha^2)^{[n/s]}\right)+O(n)\\ &=& 2\sum_{\substack{1\le s\leq n\\
            r\le s,\ (r,s)=1}}\left(1-(1-\alpha^2)^{[n/s]}\right)+O(n)\\ &=&2\sum_{s\leq n}\varphi(s)\left(1-(1-\alpha^2)^{[n/s]}\right)+O(n), 
\end{eqnarray*}     
where $\varphi$ stands for the Euler's totient function. 

Now we observe that $k\leq n/s<k+1$ if and only if $n/k\geq s> n/(k+1)$. Denoting $\Phi(x):=\sum_{n\leq x}\varphi(n)$ and summing by parts, we obtain 
\begin{align*}
 \mathbb{E}(X) 
 &= 2\sum_{k=1}^{n-1}\sum_{\frac{n}{k+1}< s\leq \frac{n}{k}}\varphi(s)\left(1-(1-\alpha^2)^{k}\right)+ O(n)\\
 &= 2\sum_{k=1}^{n-1}\left(\Phi\left(\frac{n}{k}\right)-\Phi\left(\frac{n}{k+1}\right)\right)\left(1-(1-\alpha^2)^{k}\right)+ O(n)\\
 &= 2\sum_{k=1}^{n-1}\Phi\left(\frac{n}{k}\right)\left(\left(1-(1-\alpha^2)^{k}\right)-\left(1-(1-\alpha^2)^{k-1}\right)\right)+ O(n)\\
 &= 2\alpha^2\sum_{k=1}^{n-1}\Phi\left(\frac{n}{k}\right)(1-\alpha^2)^{k-1}+ O(n).
\end{align*}
The classical estimate $\Phi(x)=\frac{3x^2}{\pi^2}+O(x\log x)$ finishes the proof:
\begin{align*} 
 \mathbb{E}(X)&= 2\alpha^2\sum_{k=1}^{n-1}\left(\frac{3n^2}{\pi^2k^2}+O\left (\frac nk\log (n/k)\right )\right)(1-\alpha^2)^{k-1}+ O(n)\\
  &= \frac{6\alpha^2n^2}{\pi^2(1-\alpha^2)}\left(\sum_{k=1}^{\infty}\frac{(1-\alpha^2)^k}{k^2}-\sum_{k=n}^{\infty}\frac{(1-\alpha^2)^k}{k^2}\right)+O(n\log^2n)\\
 &=\frac{6\alpha^2n^2}{\pi^2(1-\alpha^2)}\sum_{k=1}^{\infty}\frac{(1-\alpha^2)^k}{k^2}+O(n\log^2 n)\\
 &=\frac {6}{\pi^2}(n\alpha)^2\frac{{\rm Li}_2(1-\alpha^2)}{1-\alpha^2}+O(n\log^2n).
 \end{align*} \end{proof}        

\subsection{Variance}\label{subsec:var}
Our next step to prove Theorem \ref{TH1} is to estimate the deviation of $X$ from its expected value. To accomplish this, we obtain a bound for its variance:

\begin{prop}\label{var} For the random variable $X=|A/A|$ in $B(n,\alpha)$ we have $$\mathrm{Var}(X)\ll n^3\log^2 n.$$
\end{prop}
\begin{proof} Since $\mathrm{Var}(X)=\mathbb{E}(X^2)-\mathbb{E}^2(X)$ and
$\mathbb{E}(X)$ was explicitly computed in the previous proposition, it suffices to estimate $\mathbb{E}(X^2)$. Linearity of expectation gives the equality
\begin{equation}\label{eq3}\mathbb{E}(X^2)=\sum_{\substack{r,s,r',s'\le n\\(r,s)=(r',s')=1}}\mathbb{P}\left(\lbrace r/s\in A/A\rbrace\cap\lbrace r'/s'\in A/A\rbrace\right),\end{equation}
so we are led to study the independence of the previous events. Given $r,s,r',s'$, the events $$\lbrace r/s\in A/A\rbrace= \bigvee_{t\le n/\max(r,s)}\{rt,st\in A\}\quad \text{ and }\quad \lbrace r'/s'\in A/A\rbrace= \bigvee_{t'\le n/\max(r',s')}\{r't',s't'\in A\}$$  are dependent if and only if there are $t\leq n/\max(r,s)$ and $t'\leq n/\max(r',s')$ such that $rt=r't'$ or $rt=s't'$ or $st=r't'$ or $st=s't'$. Now $rt=r't'$ if and only if $tr/(r,r')=t'r'/(r,r')$, so coprimality implies the existence of positive integers $u,u'$ such that $t=ur'/(r,r')$ and $t'=u'r/(r,r')$. We observe that
$$\frac{n}{r}\geq t= \frac{ur'}{(r,r')}\geq \frac{r'}{(r,r')};$$
so condition $n<rr'/(r,r')$ guarantees that there are no $t\leq n/\max(r,s), t'\leq n/\max(r',s')$ such that $rt=r't'$. Similarly, conditions $n< rr'/(r,r'), sr'/(s,r'), rs'/(r,s'), ss'/(s,s')$ ensure that $\lbrace r/s\in A/A\rbrace$ and $\lbrace r'/s'\in A/A\rbrace$ are independent events. 

Consequently, we split the sum in ~\eqref{eq3} in two parts: the terms which satisfy the former conditions and the terms which do not satisfy them. In the first one, all the events are independent, so their sum is bounded by
\begin{equation}\label{eq4}
\left(\sum_{\substack{r,s\le n\\ (r,s)=1}}\mathbb{P}\left(\lbrace r/s\in A/A\rbrace\right)\right)\left(\sum_{\substack{r',s'\le n\\ (r',s')=1}}\mathbb{P}\left(\lbrace r'/s'\in A/A\rbrace\right)\right)=\mathbb{E}(X)^2.
\end{equation}
In the second one, we simply bound the probabilities by 1. Hence, we are led to count the number of irreducible fractions $r/s, r'/s'$ whose denominators and numerators are bounded by $n$ and which satisfy one of the following conditions: $rr'/(r,r'), sr'/(s,r'), rs'/(r,s'), ss'/(s,s')\leq n$. By symmetry, it suffices to treat the first case. Firstly, we observe that
$$\sum_{\substack{r, s,  r',s'\leq n\\ 
            (r,s)=(r',s')=1\\ rr'/(r,r')\leq n}} 1 \ll n^2\sum_{\substack{r, r'\leq n\\ 
            rr'/(r,r')\leq n}} 1.$$
Then, defining $l:=(r,r')$ and writing $r=lm, r'=lm'$ for certain integers $m,m'$, the expression above is bounded by
$$n^2\sum_{\substack{l\leq n\\ m, m'\leq n/l\\ 
            mm'\leq n/l}} 1\ll n^2\sum_{\substack{l\leq n\\ k\leq n/l}} \tau(k),$$           
where we denote by $\tau(k)$ the number of divisors of $k$. Finally, the classical estimates for the average order of $\tau(k)$ and for the harmonic numbers imply that the former quantity is
\begin{equation}\label{eq5}
\ll n^2\sum_{l\leq n} \frac{n}{l}\log\left(\frac{n}{l}\right)\ll n^3\log n\sum_{l\leq n} \frac{1}{l}\ll n^3\log^2 n.
\end{equation}
Substituting the estimates for ~\eqref{eq4} and ~\eqref{eq5} into ~\eqref{eq3}, we obtain
$$\mathbb{E}(X^2)\le \mathbb{E}(X)^2+O(n^3\log^2 n).$$
Thus,
$$\mathrm{Var}(X)=\mathbb{E}(X^2)-\mathbb{E}(X)^2\ll n^3\log^2 n,$$
as we wished to show.
\end{proof}

Now we are in a position to prove Theorem \ref{TH1}, since the estimate $\mathrm{Var}(X)=o(\mathbb{E}(X)^2)$ and Chebyshev's inequality imply that, for every $\epsilon>0$, 
$$\lim_{n\rightarrow\infty}{\mathbb{P}\left(\left|\frac{X}{\mathbb{E}(X)}-1\right|\geq \varepsilon\right)}= 0. $$
Hence, $X\sim\mathbb{E}(X)$ with probability $1-o(1)$ as $n\to \infty$, and Theorem \ref{TH1} is obtained.

\begin{proof}[Proof of Corollary \ref{strong}]
Using Chebyshev's inequality and Proposition \ref{var} we have
$$\sum_M\mathbb P\left (\left |\frac{|A_{M^3}/A_{M^3}|}{M^{6}}-\mathbb E\left ( \frac{|A_{M^3}/A_{M^3}|}{M^{6}} \right )\right |>1/\sqrt M \right )\le \sum_M
\frac{\mathrm{Var}\left ( \frac{|A_{M^3}/A_{M^3}|}{M^6}\right )}{(1/\sqrt M)^2}\ll \sum_M\frac{\log^2M}{M^2}.
$$
Since the former sum is convergent, the Borel-Cantelli lemma implies that
$$\left |\frac{|A_{M^3}/A_{M^3}|}{M^{6}}-\mathbb E\left ( \frac{|A_{M^3}/A_{M^3}|}{M^{6}} \right )\right |\ll 1/\sqrt M $$ almost surely.

For each positive integer $n$, let $M=M(n)$ such that $M^3\le n<(M+1)^3$. It is clear that $$|A_n/A_n|=|A_{M^3}/A_{M^3}|+O(M^5).$$ We have that, almost surely,
\begin{eqnarray*}\left |\frac{|A_n/A_n|}{n^2}-\mu_2(\alpha)\right |&\le &\left |\frac{|A_n/A_n|}{n^2}-\frac{|A_{M^3}/A_{M^3}|}{M^6}\right |+
\left |\frac{|A_{M^3}/A_{M^3}|}{M^6}-\mathbb E\left (\frac{|A_{M^3}/A_{M^3}|}{M^6}\right )\right |\\ &+& \left |\mathbb E\left (\frac{|A_{M^3}/A_{M^3}|}{M^6}\right )-\mu_2(\alpha)\right |\\ &\ll &\frac 1M+\frac 1{\sqrt M}+\left |\mathbb E\left (\frac{|A_{M^3}/A_{M^3}|}{M^6}\right )-\mu_2(\alpha)\right |\to 0\end{eqnarray*}
as $n\to \infty$. In other words, with probability $1$ we have that $\lim_{n\to \infty}\frac{|A_n/A_n|}{n^2}=\mu_2(\alpha),$ and Corollary \ref{strong} is obtained.
\end{proof}

\section{Random sets in $S(n,k)$. Proof of Theorem \ref{TH2}}
We follow the same strategy used in \cite{C}. Let us consider again the random variable $X=|A/A|$, but in the model $S(n,k)$. From now on, $\EE_k(X)$ and $V_k(X)$ will denote the expected value and the variance of $X$ in this probability space. Clearly, for $s=1,2$ we have
\begin{eqnarray*}
\EE_k(X^s)&=&\frac 1{\binom
nk}\sum_{|\A|=k }|A/A|^s\\ V_k(X)&=&\frac 1{\binom nk}\sum_{|\A|=k}\left
(|A/A|-\EE_k(X)\right
)^2\end{eqnarray*}
\begin{lemma}\label{k-j}
For $s=1, 2$  and $1\le j<k$ we have that
$$\EE_j(X^s)\le \EE_k(X^s)\le
\EE_j(X^s)+k^{2s}-j^{2s}.$$
\end{lemma}

\begin{proof} In order to prove the lower bound it is enough to consider the case $j = k-1$. Observe that $|A/A|$ is monotone with respect to inclusion, i.e. $|\left(\A\cup \{a\}\right)/\left(\A\cup \{a\}\right)|\geq |A/A|$ for any $\A,\{a\}\subseteq \{1,\ldots, n\}$.
  Using this we get
  \begin{eqnarray*}
  \mathbb E_{k-1}(X^s)&=&\frac 1{\binom n{k-1}}\sum_{|A| = k-1}|A/A|^s\\ &\leq& \frac 1{\binom n{k-1}}\sum_{|A| = k-1}\frac{1}{n-k+1}\sum_{a\in \{1,\ldots, n\}\setminus
  A}|\left(\A\cup \{a\}\right)/\left(\A\cup \{a\}\right)|^s\\ &=&\frac 1{\binom n{k-1}} \frac{k}{(n-k+1)}\sum_{|A'| = k}|\A'/\A'|^s\\ &=&\frac 1{\binom n{k}}\sum_{|A'| = k}|\A'/\A'|^s=\mathbb E_{k}(X^s).
  \end{eqnarray*}

For the second inequality, we observe that for any set $\A\subset\{1,\ldots, n\}$ of size $k$ and any partition into two sets $\A=\A'\cup \A''$ with $|\A'|=j,\
|\A''|=k-j$ we have that $|A/A|\le |A'/A'|+2|A'||A''|+|A''|^2=
|A'/A'|+k^2-j^2$. Similarly, \begin{eqnarray*}|A/A|^2&\le &
(|A'/A'|+k^2-j^2)^2\\ &= &|A'/A'|^2+2|A'/A'|(k^2-j^2)+(k^2-j^2)^2\\ &\le & |A'/A'|^2+2j^2(k^2-j^2)+(k^2-j^2)^2\\
&= &|A'/A'|^2+k^4-j^4.\end{eqnarray*} Thus, for $s=1,2$
we have
\begin{eqnarray*} |A/A|^s&\le &\binom kj^{-1} \sum_{\substack{\A'\subset \A\\
|\A'|=j}}\left
(|A'/A'|^s+k^{2s}-j^{2s}\right )\\&\le &\binom kj^{-1}
\Big(\sum_{\substack{\A'\subset \A\\ |\A'|=j}}|A'/A'|^s\Big)+k^{2s}-j^{2s}.\end{eqnarray*} Then,
\begin{eqnarray*}\sum_{|\A|=k}|A/A|^s&\le &\binom
kj^{-1}\sum_{|\A|=k}\sum_{\substack{\A'\subset \A\\
|\A'|=j}}|A'/A'|^s+\binom
nk(k^{2s}-j^{2s})\\
&= &\binom kj^{-1}\sum_{|\A'|=j}|A'/A'|^s\sum_{\substack{\A'\subset
\A\\
|\A|=k}}1+\binom nk(k^{2s}-j^{2s})\\
&= &\binom kj^{-1}\binom{n-j}{k-j}\sum_{|\A'|=j}|A'/A'|^s+\binom
nk(k^{2s}-j^{2s})\\&= &\frac{\binom nk}{\binom
nj}\sum_{|\A'|=j}|A'/A'|^s+\binom nk(k^{2s}-j^{2s}),
\end{eqnarray*}
and the second inequality holds.
\end{proof}
\begin{prop}\label{Tcomb} For $s=1,2$ we have that
$$\EE_k(X^s)=\EE(X^s)+O(k^{2s-1/2}),
$$
where $\EE(X^s)$ denotes the expectation of $X^s$  in $
B(n,k/n)$ and $\EE_k(X^s)$ the expectation in $S(n,k)$.
\end{prop}
\begin{proof}
Observe that for $s=1,2$ we have
\begin{eqnarray*}\EE(X^s)-\EE_k(X^s)&=&-
\EE_k(X^s)+\sum_{j=0}^n\left(\frac{k}{n}
\right)^j\left(1-\frac{k}{
n}\right)^{n-j}
\sum_{|\A|=j}|A/A|^s\\&=&-\EE_k(X^s)+\sum_{j=0}^n
\left(\frac{k}{n}\right)^j\left(1-\frac{k}{n}\right)^{n-j}\binom
nj\EE_j(X^s)\\
&=&\sum_{j=0}^n\left(\frac{k}{n}\right)^j\left(1-\frac{k}{n}\right)^{
n-j}\binom nj\left
(\EE_j(X^s)-\EE_k(X^s)\right),
\end{eqnarray*}
for $s=1,2$. Using Lemma \ref{k-j} we get
\begin{equation}\label{Ejk}|\EE_k(X^s)-\EE(X^s)|\le\sum_{j=0}^n\left(\frac{k}{n}\right)^j\left(1-\frac{k}{n}\right)^{
n-j}\binom nj|j^{2s}-k^{2s}|.
\end{equation}
To estimate the sum above, we observe that
\begin{equation}\label{jk}
|j^{2s}-k^{2s}|\le 4|j-k|(\max(j,k))^{2s-1}.
\end{equation}
We also consider $\mathbb{E}(|Y-\mathbb{E}(Y)|)$, where $Y\sim \text{Bin}(n,k/n)$ is the binomial distribution of parameters $n$ and $k/n$. Cauchy--Schwarz inequality for the expectation implies that this quantity is bounded by the standard deviation of the binomial distribution.
 \begin{equation}\label{ss}\sum_{j=0}^n\left(\frac{k}{n}\right)^j\left(1-\frac{k}{n}\right)^{
n-j}\binom nj|j-k|=\mathbb{E}(|Y-\mathbb{E}(Y)|)\le \sqrt{\text{Var}(Y)}= \sqrt{n(k/n)(1-k/n)}\le \sqrt k.\end{equation}

To estimate the sum in~\eqref{Ejk}  we split the expression in two terms: the sum indexed by $j\leq 2k$ and the one with $j>2k$.
We use  \eqref{jk} and \eqref{ss} to get
\begin{eqnarray*}
\sum_{j\le 2k} \left(\frac{k}{n}\right)^j\left(1-\frac{k}{n}\right)^{
n-j}\binom nj|j^{2s}-k^{2s}|&\le &4(2k)^{2s-1}\sum_{j=0}^n\left(\frac{k}{n}\right)^j\left(1-\frac{k}{n}\right)^{
n-j}\binom nj|j-k|\\ &\le & 32k^{2s-1/2}.
\end{eqnarray*}
On the other hand,
\begin{eqnarray*} & &\sum_{j> 2k}\left(\frac{k}{n}\right)^j\left(1-\frac{k}{n}\right)^{
n-j}\binom nj|j^{2s}-k^{2s}|\\ &\le & \sum_{j> 2k} \left(\frac{k}{n}\right)^j\left(1-\frac{k}{n}\right)^{
n-j}\binom njj^{2s}\\
&\le &\sum_{l\ge 2}(l+1)^{2s}k^{2s}\sum_{lk<j\le (l+1)k} \left(\frac{k}{n}\right)^j\left(1-\frac{k}{n}\right)^{
n-j}\binom nj\\
&\le  &\sum_{l\ge 2}(l+1)^{2s}k^{2s}\ \mathbb P(Y>lk)
\end{eqnarray*}
where, once again,  $Y\sim \text{Bin}(n,k/n)$. Chernoff's Theorem implies that for any $\epsilon>0$ we have
$$\mathbb P(Y>(1+\epsilon)k)\le e^{-\epsilon^2k/3}.$$ Applying this inequality to $\mathbb P(Y>lk)$ we get
\begin{eqnarray*}\sum_{j> 2k} & &\left(\frac{k}{n}\right)^j\left(1-\frac{k}{n}\right)^{
n-j}\binom nj|j^{2s}-k^{2s}|\\
\le & &\sum_{l\ge 2}(l+1)^{2s}k^{2s} e^{-(l-1)^2k/3}\ll k^{2s}e^{-k/3}\ll k^{2s-1/2}.
\end{eqnarray*}
\end{proof}
The next corollary proves the first part of Theorem \ref{TH2}.
\begin{cor}\label{coro2} If $k\asymp n$  then
$$\EE_k(X)=\frac{6}{\pi^2}k^2\frac{\mathrm{Li}_2(1-(k/n)^2)}{1-(k/n)^2}\left (1+O\left(n^{-1/2}\right) \right ).$$
\end{cor}
\begin{proof}
Proposition \ref{Tcomb} for $s=1$ and Proposition \ref{expectation} imply that 
\begin{eqnarray*}
\EE_k(X)&=&\EE(X)+O(n^{3/2})\\ &=&\frac{6}{\pi^2}k^2\frac{\text{Li}_2(1-(k/n)^2)}{1-(k/n)^2}+O(n\log^2 n)+O(k^{3/2})\\ &=&\frac{6}{\pi^2}k^2\frac{\text{Li}_2(1-(k/n)^2)}{1-(k/n)^2}\left (1+O\left(n^{-1/2}\right) \right ).
\end{eqnarray*}
\end{proof}

To conclude the proof of Theorem~\ref{TH2}, we combine Proposition \ref{var} and Proposition \ref{Tcomb} to estimate the variance $V_k(X)$ in $S(n,k)$ using the variance $V(X)$ in $B(n,k/n)$:
\begin{align*}
 V_k(X)&=\EE_k(X^2)-\EE_k^2(X)\\&=V(X)+\left
(\EE_k(X^2)-\EE(X^2)\right )+\left
(\EE(X)-\EE_k(X)\right )\left
(\EE(X)+\EE_k(X)\right )
\\ &\ll n^3\log^2n+k^{7/2}+\left (k^{3/2}\right )\left ( k^2\right )
\\ & \ll  k^{7/2}.
\end{align*}
The assertion of Theorem \ref{TH2} is a consequence of the estimate  $V_k(X)=o\left (\EE_k^2(X)\right )$ when $k\to \infty$.

\section{The size of the ratio set $A/A$ when $A$ is an arithmetic progression.}
\begin{proof}[Proof of Theorem \ref{TH3}]Without loss of generality, we may assume that $(a,q)=1$ by cancelling the factor $(a,q)$ in all the fractions of $A/A$. 

If $(r,s)=1$, note that $r/s\in A/A$ if and only if $(rt,st)\in A\times A$ for certain positive integer $t$, which occurs if and only if $rt\equiv st\equiv a\pmod{q}$ for some $t$ with $rt,st\leq n$. A necessary, but not sufficient, condition for this is that $r\equiv s\pmod{q}$. Indeed, $rt\equiv st\equiv a\pmod{q}$ and $(a,q)=1$ imply that $(t,q)=1$, so $r\equiv s\equiv t^{-1}a\pmod{q}$.  

We consequently classify $r, s$ according to their remainder $l$ modulo $q$. Thus, if $r\equiv s\equiv l \pmod{q}$ then $r/s\in A/A$  if and only if there exists $t=t(l)$ such that $lt\equiv a\pmod{q}$ and $rt,st\leq n$. For the first condition, it is enough that $(l,q)=1$; and it is necessary for, if such a $t$ exists, then $(l,q)|a$ and $(a,q)=1$. Hence
$$|A/A|=\sum_{\substack{1\leq l\leq q\\ (l,q)=1}}|S(l)|,$$ 
where $S(l):=\{r,s\in\mathbb{N} :  (r,s)=1, r\equiv s\equiv l\pmod{q},\ rt(l), st(l)\leq n \}$ and $t(l)$ is the least positive integer congruent to $l^{-1}a$ modulo $q$. 

We begin by estimating the cardinality of these sets. Observe that
\begin{equation}\label{eq100}
|S(l)|=\sum_{\substack{r, s\leq n/t(l)\\ r\equiv s\equiv l \pmod{q} \\ (r,s)=1}} 1=\sum_{\substack{r, s\leq n/t(l)\\ r\equiv s\equiv l \pmod{q} }}\sum_{d|(r,s)} \mu(d),
\end{equation}  
  and that, for $d|(r,s)$, we may write $r=dr', s=ds'$, where $r',s'\leq n/t(l)d$. Since $(q,l)=1$ and $(d,q)|l$ because $dr'=r\equiv l\pmod{q}$, we must have $(d,q)=1$. Thus, $r\equiv s\equiv l \pmod{q}$ is equivalent to $r'\equiv s'\equiv d^{-1}l \pmod{q}$, and we may rewrite ~\eqref{eq100} as
     $$|S(l)|= \sum_{\substack{d\leq n/t(l)\\ (d,q)=1}}\mu(d) \sum_{\substack{r', s'\leq n/t(l)d\\ r'\equiv s'\equiv d^{-1}l \pmod{q} }}1 =\sum_{\substack{d\leq n/t(l)\\ (d,q)=1}}\mu(d) \left(\frac{n}{t(l)dq}+O(1)\right)^2,$$     
since the number of positive integers less or equal than a certain $m$ and congruent to a certain $a$ modulo $b$ is $m/b+O(1)$. Thus,
\begin{align*}
|S(l)| &=\frac{n^2}{q^2t(l)^2}\sum_{\substack{d\leq n/t(l)\\ (d,q)=1}}\frac{\mu(d)}{d^2}+O\left(\frac{n}{qt(l)}\log(n/t(l))\right).
\end{align*}

Now we split the former sum in two parts
 $$\sum_{\substack{d\leq n/t(l)\\ (d,q)=1}}\frac{\mu(d)}{d^2}=\sum_{(d,q)=1}\frac{\mu(d)}{d^2}-\sum_{\substack{d> n/t(l)\\ (d,q)=1}}\frac{\mu(d)}{d^2}.$$                   
The second one may be bounded as follows                 
 $$\sum_{\substack{d> n/t(l)\\ (d,q)=1}}\frac{\mu(d)}{d^2}\ll \sum_{d>n/t(l)}\frac{1}{d^2}=O\left(\frac{t(l)}{n}\right).$$
 And, for the first one, note that
 $$\sum_{(d,q)=1}\frac{\mu(d)}{d^2}=\prod_{p\nmid q}\left (1-\frac 1{p^2}\right )=\frac{6}{\pi^2}\prod_{p\mid q}\left (1-\frac 1{p^2}\right )^{-1}.$$ 
 Hence
$$|S(l)|=\frac{n^2}{q^2t(l)^2}\frac{6}{\pi^2}\prod_{p\mid q}\left (1-\frac 1{p^2}\right )^{-1}+O\left(\frac{n\log n}q\right ).$$
Finally, we have to add the previous quantities for all $1\leq l\leq q$ such that $(l,q)=1$. But $t(l)$ was the least positive integer congruent to $l^{-1}a$ modulo $q$, with $(a,q)=1$, and $\psi: l\mapsto l^{-1}a$ is a set automorphism in $(\mathbb{Z}/q\mathbb{Z})^{\times}$. Therefore, $\{l\in\mathbb{Z} : 1\leq l\leq q, (q,l)=1\}$ and  $\{t(l): 1\leq l\leq q, (q,l)=1\}$ coincide.  Thus,
$$|A/A|=\sum_{\substack{l\leq q\\ (l,q)=1}}S(l)=\frac{n^2}{q^2}\frac{6}{\pi^2}\prod_{p\mid q}\left (1-\frac 1{p^2}\right )^{-1}\sum_{\scriptstyle 1\le t\le q\atop\scriptstyle (t,q)=1}\frac 1{t^2}+O(n\log n),$$
and  Theorem \ref{TH3} is proved.
\end{proof}
To prove that $c_q\to 1$ when $q\to \infty$, observe that
\begin{eqnarray*}\sum_{\scriptstyle 1\le t\le q\atop\scriptstyle (t,q)=1}\frac 1{t^2}&=&\sum_{t=1}^q\sum_{d|(t,q)}\frac{\mu(d)}{t^2}=\sum_{d|q}\frac{\mu(d)}{d^2}\sum_{t'\leq q/d}\frac{1}{t'^2}=\sum_{d|q}\frac{\mu(d)}{d^2}\left(\frac{\pi^2}{6}+O\left(\frac{d}{q}\right)\right)\\&=&\frac{\pi^2}6\prod_{p\mid q}\left (1-\frac 1{p^2}\right )+O\left (\frac 1q\sum_{d\mid q}\frac 1d\right ).\end{eqnarray*}

We note that $\frac {1}{q}\sum_{d|q}\frac 1d\ll \frac {\tau(q)}{q}$ where $\tau(q)$ is the number of divisors of $q$, so
$$c_q=\frac{6}{\pi^2}\prod_{p\mid q}\left (1-\frac 1{p^2}\right )^{-1}\sum_{\scriptstyle 1\le t\le q\atop\scriptstyle (t,q)=1}\frac 1{t^2}=1+O\left(\frac {\tau(q)}{q}\right)\xrightarrow{\: q \to \infty \: } 1,$$
since $\frac{6}{\pi^2}\prod_{p\mid q}\left (1-\frac 1{p^2}\right )^{-1}=\prod_{p \nmid q}\left (1-\frac 1{p^2}\right )\le 1$.

\section{Bounds for $|A/A|$ when $A$ has positive density}
We now pass to the proof of Theorem \ref{TH5}.
\begin{proof}
Given $\epsilon>0$, let $m$ be a positive integer such that $(3/4)^m<\epsilon$. We shall denote by $\mathcal{P}=\{p_1,\dots,p_m\}$ the set of the first $m$ prime numbers, and by $S_m=\{p_1^{\epsilon_1}\cdots p_m^{\epsilon_m}:\ \epsilon_i\in \{0,1\}\}$ the set of all the products of different elements of $\mathcal{P}$. Furthermore, given a large integer $n>m$, we consider the set $$A':=\{r\in\mathbb{N}: r\leq n/(p_1\cdots p_m), (r,\mathcal{P})=1\}$$ and the set $$A:=\bigcup_{s\in S_m} sA'=\{sr: s\in S_m, r\leq n/(p_1\cdots p_m), (r,\mathcal{P})=1\}.$$
Since this is a disjoint union, we have that 
\begin{equation}\label{disjoint}
|A|=|S_m||A'|.
\end{equation}
Moreover,
$$A/A=\left \{\frac s{s'}\frac r{r'}:\ s/s'\in S_m/S_m,\ r/r'\in A'/A',\ (s,s')=(r,r')=1\right \};$$ and, since $(s,s')=(r,r')=1$ and $(r,\mathcal{P})=(r',\mathcal{P})=1$, all the products  $\frac s{s'}\frac r{r'}$ in $A/A$ are distinct. Hence
\begin{equation}\label{dis2}
|A/A|=|S_m/S_m||A'/A'|.
\end{equation}
Putting \eqref{disjoint} and \eqref{dis2} together, we get
$$|A/A|=|S_m/S_m||A'/A'|\le |S_m/S_m||A'|^2 = \frac{|S_m/S_m|}{|S_m|^2}|A|^2.$$
Now $|S_m/S_m|=|\{p_1^{\epsilon_1}\cdots p_m^{\epsilon_m}:\ \epsilon_i\in \{-1,0,1\}\}|=3^m$ and $|S_m|=2^m$, so $\frac{|S_m/S_m|}{|S_m|^2}=(3/4)^m<\epsilon $. Thus,
$$|A/A|<\epsilon |A|^2.$$

Finally, using that 
$$|\{ r\le x:\ (r,M)=1\}|=\frac{\varphi(M)}Mx+O(\tau(M)),$$ we obtain
$$|A|=|S_m||A'|=2^m\left (\frac n{p_1\cdots p_m}\frac{\varphi(p_1\cdots p_m)}{p_1\cdots p_m}+O(2^m)\right )\ge \alpha n(1+o(1)),$$
where $\alpha:=\frac{2^m\varphi(p_1\cdots p_m)}{(p_1\cdots p_m)^2}$>0. 

This concludes the proof of the theorem and leads us into the last section of this paper.
\end{proof}

\section{Visible lattice points in multidimensional spaces. Proof of Theorem \ref{TH4}}

In order to prove Theorem \ref{TH4}, it is convenient to introduce a certain generalization of the Euler's totient function.

\begin{defi} Given $k, m\in\mathbb{N}$, the Jordan's totient function $J_k(m)$  is defined as the cardinality of the set $\lbrace(a_1,...,a_k)\in\mathbb{N}^k : a_1\leq...\leq a_k\leq m, \mathrm{gcd}(a_1,...,a_k,m)=1\rbrace$. 
\end{defi}

Plainly, $J_1=\varphi$. Further, it is known (\cite{SS}) that
\begin{equation}\label{Jordan}
\Phi_k(n):=\sum_{m\le n}J_k(m)=\frac{n^{k+1}}{(k+1)!\zeta(k+1)}+O(n^k).
\end{equation}

The proof of Theorem \ref{TH4} is now similar to that of Theorem \ref{TH1}, so we shall only indicate the appropriate modifications.

\begin{proof}
Let $X$ denote the random variable $|{\rm visible}(A_1\times \cdots \times A_d)|$, where each element of $A_i$ is chosen independently at random in $\{1,\dots, n\}$ with probability $\alpha_i$, for $1\leq i\leq d$. By symmetry,
\begin{align}\label{eq50}
\mathbb{E}(X) &=\sum_{\substack{a_1,...,a_d\leq n\\
            (a_1,...,a_d)=1}}\mathbb{P}\left(\lbrace (a_1t,...,a_dt)\in A_1\times \cdots \times A_d\  \text{for some } t\leq n/\max(a_i)\rbrace\right) \notag\\
            &=d!\sum_{\substack{a_1<...<a_d\leq n\\
            (a_1,...,a_d)=1}}\mathbb{P}\left(\lbrace (a_1t,...,a_dt)\in A_1\times \cdots \times A_d\  \text{for some } t\leq n/a_d)\rbrace\right) +O(n^{d-1}),
\end{align}   
where $O(n^{d-1})$ bounds the number of points having at least two identical coordinates. If $a_1<\cdots <a_d$, we have
\begin{eqnarray*}& &\mathbb{P}\left(\lbrace (a_1t,...,a_dt)\in A_1\times \cdots \times A_d\  \text{for some } t\leq n/a_d\rbrace\right)\\ &=&1-\mathbb{P}\left(\lbrace (a_1t,...,a_dt)\not \in A_1\times \cdots \times A_d\  \text{for any } t\leq n/a_d)\rbrace\right)\\&=&1-\prod_{t\le n/ a_d}\mathbb{P}\left(\lbrace (a_1t,...,a_dt)\not \in A_1\times \cdots \times A_d\rbrace\right)\\&=&1-\prod_{t\le n/a_d}\left (1-\mathbb{P}\left(\lbrace (a_1t,...,a_dt)\in A_1\times \cdots \times A_d\rbrace\right)\right )\\ &=&1-\left (1-\alpha_1\cdots\alpha_d\right )^{[ n/a_d]}.\end{eqnarray*}

Thus, we may rewrite ~\eqref{eq50} as
\begin{align*}
\mathbb{E}(X) &=d!\sum_{\substack{a_1<\dots <a_d\leq n\\
            (a_1,...,a_d)=1}}\left(1-(1-\alpha_1\cdots \alpha_d)^{[n/a_d]}\right)+O(n^{d-1})\\
            &=d!\sum_{a_d\leq n}J_{d-1}(a_d)\left(1-(1-\alpha_1\cdots \alpha_d)^{[n/a_d]}\right )+O(n^{d-1}).
            \end{align*}            
            
            Finally, the explicit expression for $\mathbb{E}(X)$ follows by proceeding as in the last part of the proof of Proposition 2.1. and by applying \eqref{Jordan}:
       \begin{align*}
\mathbb{E}(X) &=d!\sum_{k=1}^{n-1}\sum_{\frac n{k+1}<a_d\leq \frac nk}J_{d-1}(a_d)\left(1-(1-\alpha_1\cdots \alpha_d)^{k}\right )+O(n^{d-1})\\
&=d!\sum_{k=1}^{n-1}\left ( \Phi_{d-1}\left(\frac nk \right )-\Phi_{d-1}\left (\frac n{k+1}\right )  \right )\left(1-(1-\alpha_1\cdots \alpha_d)^{k}\right )+O(n^{d-1})\\&=d!\sum_{k=1}^{n-1} \Phi_{d-1}\left(\frac nk \right )\left (\left(1-(1-\alpha_1\cdots \alpha_d)^{k}\right )-\left(1-(1-\alpha_1\cdots \alpha_d)^{k-1}\right )\right ) +O(n^{d-1})\\&= d!\alpha_1\cdots\alpha_d\sum_{k=1}^{n-1} \Phi_{d-1}\left(\frac nk \right )\left(1-\alpha_1\cdots \alpha_d\right )^{k-1} +O(n^{d-1})\\&=
d!\alpha_1\cdots\alpha_d\sum_{k=1}^{n-1} \left (\frac{n^d}{d!\zeta(d)k^d}+O((n/k)^{d-1}) \right )\left(1-\alpha_1\cdots \alpha_d\right )^{k-1} +O(n^{d-1})\\&=n^d\frac{\alpha_1\cdots\alpha_d}{\zeta(d)}\sum_{k=1}^{n-1}\frac{\left(1-\alpha_1\cdots \alpha_d\right )^{k-1}}{k^d}+O(n^{d-1}\log n),
            \end{align*}      
            
             obtaining 
            $$\mathbb{E}(X)=n^d\frac{\alpha_1\cdots \alpha_d}{\zeta(d)}\frac{{\rm Li}_d(1-\alpha_1\cdots \alpha_d)}{1-\alpha_1\cdots \alpha_d}+O(n^{d-1}\log n) .$$
           
To estimate the variance of $X=|\text{visible}(A\times \cdots A_d)                      |$, we compute \begin{eqnarray*}\mathbb E(X^2)& =&\sum_{\substack{a_1,\dots,a_d,a_1',\dots,a_d'\le n\\ (a_1,\dots,a_d)=(a_1',\dots,a_d')=1}}\mathbb P\left (\begin{matrix}(ta_1,\dots,ta_d)\in A_1\times \cdots \times A_d\text{ and }(t'a_1',\dots,t'a_d')\in A_1\times \cdots \times A_d)\\ \text{ for some } t,t'\end{matrix}\right )\\ &=&  (d!)^2\sum_{\substack{a_1<\cdots <a_d\le n\\ a_1'<\cdots <a_d'\le n\\ (a_1,\dots,a_d)\\=(a_1',\dots,a_d')=1}}\mathbb P\left (\begin{matrix}&(ta_1,\dots,ta_d)\in A_1\times \cdots \times A_d,\text{ for some }t\le \frac n{a_d}\\&(t'a_1',\dots,t'a_d')\in A_1\times \cdots \times A_d,\text{ for some } t'\le \frac n{a_d'}\end{matrix}\right )+O(n^{2d-1}) \\&=& (d!)^2\sum_{\substack{a_1<\cdots <a_d\le n\\ a_1'<\cdots <a_d'\le n\\ (a_1,\dots,a_d)\\=(a_1',\dots,a_d')=1}}\mathbb P\left (ta_i,t'a_i'\in A_i,\ i=1,\dots ,d,\  \text{ for some } t\le \frac n{a_d},t'\le \frac n
 {a
 _d'}\right )+O(n^{2d-1}) .\end{eqnarray*}

As we did in the proof of Theorem \ref{TH1}, we can check that the condition $n<a_ia_i'/(a_i,a_i')$ implies that there are no $t\le n/a_d,t'\le n/a_d'$ such that $ta_i=t'a_i'$. Hence, we split the sum above in two parts. In the first one, we sum over all $a_1,\dots,a_d,a_1',\dots,a_d'\le n$ with $(a_1,\dots,a_d)=(a_1',\dots,a_d')=1$ and $a_ia_i'/(a_ia_i')>n$; and, by independence, this part is bounded by
\begin{eqnarray*} (d!)^2& &\sum_{\substack{a_1<\cdots <a_d\le n\\ (a_1,\dots,a_d)=1}} \mathbb P\left (ta_i\in A_i,\ i=1,\dots ,d,\  \text{ for some } t\le \frac n{a_d}\right )\\  &\times &\sum_{\substack{a_1'<\cdots <a_d'\le n\\ (a_1',\dots,a_d')=1}}\mathbb P\left (t'a_i'\in A_i,\ i=1,\dots ,d,\  \text{ for some } t'\le \frac n{a_d'}\right )+O(n^{2d-1})\\&=&\mathbb E^2(X)+O(n^{2d-1}) .\end{eqnarray*}
In the second one, we sum over all $a_1,\dots,a_d,a_1',\dots,a_d'\le n$ with $a_ia_i'/(a_i,a_i')\le n$ for some $i=1,\dots,n$. This is clearly bounded by
$$d\sum_{\substack{a_1,\dots,a_d\le n\\a_1',\dots,a_d'\le n\\ a_1a_1'/(a_1a_1')\le n}}1\le dn^{2d-2}\sum_{\substack{a_1,a_1'\le n\\a_1a_1'/(a_1,a_1')\le n }}1.$$
On the other hand,
$$\sum_{\substack{a_1,a_1'\le n\\a_1a_1'/(a_1,a_1')\le n }}1=\sum_{l\le n} \sum_{\substack{a_1,a_1'\le n\\ (a_1,a_1')=l\\a_1a_1\le ln }}1\le  \sum_{l\le n} \sum_{\substack{b_1,b_1'\le n/l\\b_1b_1'\le n/l }}1\le\sum_{l\le n} \sum_{\substack{m\le n/l }}\tau(m)\ll \sum_{l\le n}\frac nl \log(n/l)\ll n\log^2n.$$
Thus, the second sum is $O\left (n^{2d-1}\log^2n\right )$, so we have that
$$\mathrm{Var}(X)=\mathbb E(X^2)-\mathbb E^2(X)\ll n^{2d-1}\log^2n.$$
Since $\mathrm{Var}(X)=o(\mathbb E^2(X))$, then 
$X\sim \mathbb E(X)$ with probability $1-o(1)$ as $n\to \infty$ and the proof is concluded.
\end{proof}

\subsubsection*{Acknowledgments:} This work was supported by grants MTM
2014-56350-P of MINECO and ICMAT Severo Ochoa project SEV-2011-0087.

\end{document}